\def\draft#1{}                      

\documentclass[11pt]{article}

\newtheorem{theorem}{Theorem}
\newtheorem{lemma}[theorem]{Lemma}

\newtheorem{cor}[theorem]{Corollary}
\newtheorem{definition}[theorem]{Definition}
\newtheorem{conjecture}[theorem]{Conjecture}

\long\def\rests#1{}

\def\noi{\noindent}

\def\pf{\noi{\bf Proof.\ \,}}
\def\eop{\hfill\framebox[2.4mm][t1]{\phantom{x}} \vskip 0.15cm } 

\def\sss#1{\if#1..\ \else\if#1,,\ \else\ #1\fi\fi} 

\def\voa#1{vertex operator algebra\sss{#1}}
\def\voas#1{vertex operator algebras\sss{#1}}
\def\svoa#1{vertex operator superalgebra\sss{#1}}
\def\svoas#1{vertex operator superalgebras\sss{#1}}
\def\subvoa#1{vertex operator subalgebra\sss{#1}}

\def\C{{\bf C}}

\def\Z{{\bf Z}}

\def\tr{{\rm tr}}

\def\halb{\frac{1}{2}\phantom{\frac{|}{|}}\!\!\!}
\def\ha{\frac{1}{2}\phantom{\frac{|}{|}}\!\!\!}
\def\bis{\hbox{--}}
\def\N1{$N\!\!=\!\!1$}

\title{Self-Dual Vertex Operator Superalgebras \\ of Large Minimal Weight}

\author{Gerald~H\"ohn\thanks{Department of Mathematics, 
Kansas State University, 138 Cardwell Hall, Manhattan, KS 66506-2602, USA.
 E-mail: {\tt gerald@math.ksu.edu}
}}

\date{}

\begin{document}

\maketitle

\abstract{The new general upper bound  $\mu\leq \left[\frac{c}{24}\right] + 1$ for 
the minimal weight $\mu$ of a self-dual \svoa of central 
charge $c\not=23\ha$ is proven. For central charges $c\leq 48$, further improved estimates
are given and examples of \svoas with large minimal weight are discussed.
We also study the case of \svoas with
\N1 supersymmetry which was first considered by Witten in connection
with three-dimensional quantum gravity. The upper bound  $\mu^*\leq \halb\left[\frac{c}{12}\right]+\halb$
for the minimal superconformal weight is obtained for $c\not=23\ha$.

In addition, we show that it is impossible that
the monster sporadic group acts on an extremal self-dual \N1 supersymmetric \svoa
of central charge~$48$ in a way proposed by Witten if certain standard assumptions 
about orbifold constructions hold. The same statement holds for extremal self-dual \voas of central charge~$48$.}


\section{Introduction}

Extremal self-dual \voas and superalgebras have been defined in \cite{Ho-dr}, Chapter~5.
Extremal refers here to the property that the degree of a Virasoro highest 
weight vector different from the vacuum vector must be larger then certain bounds obtained
from conditions on the characters.
The smallest such degree is called the minimal weight.
For small values of the central charge~$c$ several examples with interesting automorphism 
are known, like the Moonshine module $V^{\natural}$~\cite{FLM} for $c=24$ and the shorter
Moonshine module  $V\!B^{\natural}$ for $c=23\halb$~\cite{Ho-dr}. The notion extremal is analogous
to similar ones for binary codes and for lattices. These two cases have been studied much more intensively
because of their relations with more geometric problems and their applications to data
processing and transmissions, cf., for example,~\cite{Stone-space}. Examples
of such codes and lattices are known for lengths respectively dimensions up to about~$100$.

Recently it was shown by the author that extremal \voas can be used to construct conformal 
$t$-designs~\cite{Ho-conformal}, an algebraic structure sharing many properties with 
classical block designs and spherical designs. In another development, Witten proposed
that extremal \voas can be used to describe three-dimensional quantum gravity with
negative cosmological constant~\cite{Witten-3dgravity}. He also considers supergravity.
Extremal \voas have been further investigated in~\cite{Mansch-ext,GaiYin-ext,Gaber-ext,greek-ext}.

\medskip

The present paper continuous the study of self-dual \svoas of large minimal weight initiated in~\cite{Ho-dr}.
In Section~\ref{improved}, the previous upper bound $\mu\leq \halb \left[\frac{c}{8}\right] +\halb$
(\cite{Ho-dr}, Cor.~5.3.3)  for the minimal weight $\mu$ of a self-dual \svoa is improved in 
Theorem~\ref{newbound} to $$\mu\leq \left[\frac{c}{24}\right] + 1 $$
for central charges $c\not =23\halb$. This upper bound can sometimes be improved further;
Table~\ref{minsvoas} lists our results for $c\leq 48$.  For $c\leq 24$, there are always examples meeting 
the bound. For $24\ha\leq c\leq 48$, we discuss examples which yield the given lower bounds.

In Section~\ref{N1svoas}, we study the minimal superconformal weight of a self-dual \N1 supersymmetric \svoa.
We define the minimal superconformal weight~$\mu^*$ as the smallest positive degree of a highest weight vector 
for the \N1 super Virasoro algebra. In Theorem~\ref{n1upperbound} we obtain for $c\not=23\ha$ the upper bound
$$\mu^*\leq \halb\left[\frac{c}{12}\right]+\halb,$$
which for central charges divisible by~$12$ was found in~\cite{Witten-3dgravity}.
A self-dual \N1 supersymmetric \svoa with minimal superconformal weight meeting this bound is called extremal.

In the final section, we take a closer look on the case of central charge~$c=48$.
The relation between self-dual \svoas of minimal weight~$5/2$ and extremal
self-dual \N1 supersymmetric \svoas with extremal self-dual \voas is studied. This
allows us for $c=48$ to conclude that under reasonable assumptions it is impossible for the 
monster simple group to act by automorphisms on an extremal self-dual vertex operator algebra or 
an extremal self-dual \N1 supersymmetric superalgebra.

\smallskip 

It remains an open problem if  self-dual \voas or (\N1 supersymmetric) \svoas with minimal
(superconformal) weight larger than~$2$ exist.

\bigskip

In the rest of the introduction, we give precise definitions and discuss
the required results about \svoas. We assume that the reader is familiar with the
general notation of \svoas.

\smallskip 

The Virasoro algebra is the complex Lie algebra spanned by $L_n$, $n\in {\bf Z}$, 
and the central element $C$ with Lie bracket
\begin{equation}\label{virasoro} 
[L_m,L_n]=(m-n)L_{m+n}+\frac{m^3-m}{12}\delta_{m+n,0}\,C
\end{equation}
where $\delta_{k,0}=1$ if $k=0$ and $\delta_{k,0}=0$ otherwise.
For a pair $(c,h)$ of complex numbers the Verma module $M(c,h)$ is 
a representation of the Virasoro algebra generated by 
a highest weight vector $v\in M(c,h)$ with $Cv=cv$, $L_0v=hv$ and $L_nv=0$ for $n\geq 1$.
For $h=0$, the module $M(c,0)$ has a quotient isomorphic to $M(c,0)/M(c,1)$.

We {\it assume\/} that the \voas $V$ in this paper are
isomorphic to a direct sum of highest weight modules for the Virasoro algebra,
i.e., one has $$V=\bigoplus_{i\in I } M_i,$$
where each $M_i$ is a quotient of a Verma modules $M(c,h)$ with $h\in \Z_{\geq 0}$.
One has therefore a natural decomposition   
\begin{equation}\label{natdecomp}
V= \bigoplus_{h=0}^{\infty}\overline{M}(h)
\end{equation}
where $\overline{M}(h)$ is a direct sum of finitely many quotients of the Verma module $M(c,h)$.
The module $\overline{M}(0)$ is the \subvoa of $V$ generated by $\omega$ which we denote
also by $V_{\omega}$ and is  a quotient of $M(c,0)/M(c,1)$.
The smallest $h>0$ for which $\overline{M}(h)\not= 0$ was called in~\cite{Ho-dr}
the {\it minimal weight\/} of $V$ and denoted by $\mu(V)$. (If no such $h>0$ exists, we let
$\mu(V)=\infty$.)

\smallskip
A \voa is called {\it rational\/} (cf.~\cite{DLM-twistedrep}) if every admissible module 
is completely reducible. In this case
there are only finitely many irreducible admissible modules up to isomorphism and
every irreducible admissible module is an ordinary module. A \voa is called {\it simple\/} 
if it is irreducible as a module over itself. 

For an irreducible module $W$ there exists an $h$ such that $W=\bigoplus_{k\in \Z_{\geq 0}} W_{k+h}$
with $W_h\not= 0$, where  the degree $n$ subspace $W_n$ is the eigenspace of $L_0$ for the eigenvalue~$n$.
We call  $h$ the {\it conformal weight\/} of the module $W$.

\smallskip

The {\it character\/} of a module $W$ of conformal weight~$h$ 
is defined by 
$$\chi_W=q^{-c/24}\sum_{k\in \Z_{\geq 0}}\dim {W_{k+h}}q^{k+h}.$$
If $V$ is assumed to be rational and satisfying the {\it $C_2$-cofiniteness condition\/}\linebreak
$\dim(V/\hbox{Span}\{x_{(-2)}y\mid x,y\in V\})<\infty$ it is a result of Zhu~\cite{Zhu-dr} that
$\chi_W$ is a holomorphic function on the complex upper half plane
in the variable $\tau$ for $q=e^{2\pi i \tau}$. We assume in this paper that the
$C_2$-cofiniteness condition is satisfied.
The family $\{\chi_W\}_W$, where $W$ runs through the isomorphism classes of
irreducible $V$-modules $W$, transforms as a vector-valued modular function
for the modular group ${\rm SL}_2(\Z)$ acting on the upper half plane in the usual
way.

\medskip

A rational \voa $V$ is called {\it self-dual\/} (other authors use the notation holomorphic
or meromorphic) if the only irreducible $V$-module 
is $V$ itself. It follows from the above mentioned result of Zhu that
the character~$\chi_V$ is a weighted homogeneous polynomial of weight~$c$
in $\sqrt[3]{j}$ (given the weight~$8$) and $1$ (weight~$24$) where $\sqrt[3]{j}$ is the third 
root of the elliptic modular function $j$ (cf.~\cite{Ho-dr}, Thm.~2.1.2).
In particular, the central charge~$c$ of a self-dual \voa is divisible by~$8$. 
It was shown in~\cite{Ho-dr}, Cor.~5.2.3, that the minimal weight of a self-dual \voa
satisfies $\mu(V)\leq [c/24]+1$. A self-dual \voa meeting this bound is called {\it extremal.}

We also need an {\it unitary\/} condition. Sufficient for this paper is to assume that
that $V$ has a real form with positive-definite invariant bilinear form and the conformal weights
of all irreducible $V$-modules are nonnegative.

\medskip

We call a self-dual \svoa $V=V_{(0)}\oplus V_{(1)}$ rational if its even \subvoa $V_{(0)}$ is rational 
and has the same associated modular braided tensor category as the even \subvoa of
$V_{\rm Fermi}^{\otimes c}$, where $V_{\rm Fermi}$ is the \svoa of central charge 
$\halb$ describing a single fermion. One has $c\in \halb \Z$ (see~\cite{Ho-dr}, Thm.~2.2.2). The
fusion algebra of $V_{(0)}$ is for $c\in2\Z$ isomorphic to $\Z[\Z_2\times \Z_2]$,
for $c\in \Z\setminus 2\Z$ isomorphic to $\Z[\Z_4]$ and for $c\in \halb\Z\setminus \Z$ 
isomorphic to the Ising fusion algebra (\cite{Ho-dr}, Thm.~2.2.5). 
The three ($c$ non-integral) or four ($c$ integral) types of irreducible $V_{(0)}$-modules 
$V_{(0)}$, $V_{(1)}$, $V_{(2)}$ (and $V_{(3)})$
have the conformal weights $0$ (for $V_{(0)}$),  $\halb\!\!\!\pmod{1}$ (for $V_{(1)}$) and $c/8\!\!\!\pmod{1}$
(for $V_{(2)}$ and, in case of four modules, $V_{(3)}$). The $V$-module $V_{(2)}\oplus V_{(3)}$ 
for integral~$c$, respectively $V_{(2)}$ for $c$ nonintegral, is called the {\it shadow\/} of $V$
and denoted by $V'$.

We call two self-dual \voas $W$ and $\widetilde{W}$ neighbours if there exists a rational self-dual
\svoa $V=V_{(0)}\oplus V_{(1)}$ such that $W=V_{(0)}\oplus V_{(2)}$ and $W=V_{(0)}\oplus V_{(3)}$
where $V_{(0)}$, $V_{(1)}$, $V_{(3)}$ and $V_{(3)}$ are the four irreducible $V_{(0)}$-modules as
above. The pairs $V$ and $W$ as well as $V$ and $\widetilde{W}$ are also called neighbours. For central
charge $c$ divisible by $8$ a self-dual \svoa $V=V_{(0)}\oplus V_{(1)}$ has the two neighbour 
\voas $V_{(0)}\oplus V_{(2)}$ and  $V_{(0)}\oplus V_{(3)}$ which could be isomorphic.

A \svoa $V_{(0)}\oplus V_{(1)}$ comes with a natural involutive automorphism $\sigma$ which acts by 
$+1$ on $V_{(0)}$ and $-1$ on $V_{(1)}$.

The character of a self-dual rational unitary \svoa has the form
\begin{equation}\label{svoachar}
\chi_V=\sum_{r=0}^ka_r\chi_{1/2}^{2c-24r}
\end{equation}
where
\begin{equation}\label{chi12}
\chi_{1/2}=q^{-1/48}\prod_{n=0}^{\infty}(1+q^{n+1/2}),
\end{equation}
$k=\left[\frac{c}{8}\right]$ and  the $a_0$, $\ldots$, $a_k$ are uniquely determined integers
(\cite{Ho-dr}, Thm.~2.2.3). 
The character of the shadow is
\begin{equation}\label{shadowachar}
\chi_{V'}=\alpha\, \sum_{r=0}^ka_r\widetilde{\chi}_{1/2}^{2c-24r}
\end{equation}
with $\alpha=1$ for integral  and $\alpha=1/\sqrt{2}$ for nonintegral values of $c$
and 
\begin{equation}\label{chis12}
\widetilde{\chi}_{1/2}=\sqrt{2}q^{1/24}\prod_{n=0}^{\infty}(1+q^n).
\end{equation}
If we let $q=e^{2\pi i\tau}$ with $\tau$ in the complex upper half-pane, then $\chi_{1/2}$, 
and hence $\chi_V$, can be considered as the Fourier expansion  of a modular function for the 
modular group $\Gamma_{\theta}=\langle S,\, T^2\rangle$ in the cusp $i\infty$
where $S=\left( 0 \,-1 \atop 1 \,\phantom{-}0\right)$ and $T=\Bigl({{1\, 1} \atop { 1\, 0}}\Bigr)$.
The character $\chi_{V'}$ is then the expansion of $e^{2\pi ic/24}\chi_{V}$
($c$ integral), respectively  $e^{2\pi ic/24}\chi_{V}/\sqrt{2}$ ($c$~non-integral),
in the other cusp of $\overline{{\bf H}/\Gamma_{\theta}}$ represented by~$1$.


\section{Self-dual \svoas{}}\label{improved}

First we prove our new general upper bound for the minimal weight of a self-dual \svoa.
\begin{theorem}\label{newbound}
A self-dual \svoa $V$ of central charge $c$ has minimal weight
$$\mu(V)\leq \left[\frac{c}{24}\right] + 1 $$
unless $c=23\ha$ in which case $\mu(V)\leq 3/2$.
\end{theorem}
\pf
We can assume that $V$ is not a \voa, i.e.,~$V_{(1)}\not=0$, because
for \voas the estimate for $\mu(V)$ given in the theorem
was obtained in~\cite{Ho-dr}, Corollary~5.2.3.

\medskip
As discussed in the introduction, the character of $V$ is a Laurent polynomial in $\chi_{1/2}$
of the form
\begin{equation}\label{VOAentw}
\chi_V=\sum_{r=0}^ka_r\chi_{1/2}^{2c-24r}=q^{-c/24}\sum_{n=0}^{\infty} C_{n}\,q^{n/2}
\end{equation}
The character of the shadow equals
\begin{equation}\label{schattenentw}
\chi_{V'}=\alpha\, \sum_{r=0}^ka_r\widetilde{\chi}_{1/2}^{2c-24r}=
    \alpha\,  q^{c/12-[c/8]}\sum_{n=0}^{\infty} B_n\,q^n, 
\end{equation}
with $\alpha=1$ for integral  and $\alpha=1/\sqrt{2}$ for nonintegral values of $c$.

\medskip

First we consider the case $c<32$. In the range $\halb\leq c<24$, we have to show that 
$\mu(V)\leq 1$ unless $c=23\ha$. From equation~(\ref{VOAentw}) one sees directly that for 
$\halb\leq c<8$ one has $\mu(V)=\halb$. For $8\leq c<16$, the condition $C_0=1$
and $C_{1}=0$ determines $C_2=\dim V_1$ which turns out to be positive (cf.~Table~5.3 and
Table~5.4 of~\cite{Ho-dr}). Hence $\mu(V)=1$ as all vectors in $V_1$ are Virasoro highest
weight vectors. For $16\leq c<23\ha$, the condition $C_0=1$
and $C_{1}=C_2=0$ determines $\chi_V$ and thus $\chi_{V'}$.  
Unless $c=23\ha$, the character $\chi_{V'}$ has nonintegral coefficients (cf.~Table~5.4 of~\cite{Ho-dr})
and hence $\mu(V)>1$ is impossible. For $c=23\ha$ one gets $C_{3}=\dim V_{3/2} =4371$
and hence $\mu(V)=3/2$ as again all vectors in $V_{3/2}$ are Virasoro highest weight vectors.
Finally, for $24\leq c<31\ha$, we have to show that $\mu(V)\leq 2$. The condition $C_0=1$
and $C_{1}=C_2=C_{3}=0$ determines $C_4=\dim V_2$. In all cases one has $C_4>1$ and
hence $\mu(V)= 2$ as $C_4-1$ is the dimension of the space of Virasoro  highest weight vectors.

\medskip

For $c\geq 32$, the proof will only use that the coefficients $B_i$ of the characters of 
the shadow $V'$ are nonnegative rational numbers.
Let $m=\left[\frac{c}{24}\right] + 1$ and suppose that $\mu(V)>m$. 
We denote by $\chi_{M_c}=q^{-c/24}\prod_{n=2}^{\infty}(1-q^n)^{-1}$ the character of
$V_{\omega}=M(c,0)/M(c,1)$. 
Then
$$\chi_V=\chi_{M_c}\cdot (1+A_{2m+1}\,q^{m+1/2}+A_{2m+2}\,q^{m+1}+\cdots)$$
which determines the  $a_i$ for $0\leq i\leq 2m$. We will show that $a_{2m}<0$. On the
other hand, equation~(\ref{schattenentw}) allows us to write $a_{2m}$ as a linear
combination of the $B_i$ for $0\leq i \leq k-2m$, say  $a_{2m}=\sum_{i=0}^{k-2m}\beta_i B_i$.
We will show that the $\beta_i$ are all nonnegative, and thus $a_{2m}\geq0$, a contradiction.
Hence the assumption $\mu(V)>m$ must be wrong and $\mu(V)\leq \left[\frac{c}{24}\right] + 1$.

\smallskip
To determine $a_{2m}$, we let $p=q^{1/2}$ and expand $\chi_{M_c}\cdot \chi_{1/2}^{-2c}$ in
powers of $\phi=\chi_{1/2}^{-24}=p+O(p^2)$. We get
$$
\chi_{M_c}\cdot \chi_{1/2}^{-2c}=\sum_{r=0}^{\infty}\alpha_r\phi^r,
$$
where by the B\"urmann-Lagrange theorem the coefficient $\alpha_r$ is given for $r>0$ by the coefficient
of $p^{r-1}$ in
\begin{equation}\label{ar}
\frac{1}{r}\frac{d(\chi_{M_c}\cdot \chi_{1/2}^{-2c})}{dp}\left(\frac{p}{\phi}\right)^r=
\frac{1}{r}\,p^r\cdot\chi_{1/2}^{24r-2c-1}\left[\chi'_{M_c}\chi_{1/2}-2c\cdot \chi_{M_c}\chi'_{1/2}\right]
\end{equation}
and $\alpha_0=1$. Since
$$\sum_{r=0}^ka_r\phi^r=\chi_V\cdot \chi_{M_c}^{-1}\cdot \chi_{M_c}\cdot \chi_{1/2}^{-2c}=
\left(1+\sum_{n=2m+1}^{\infty} A_np^n\right)\left(\sum_{r=0}^{2m}\alpha_r\phi^r+\sum_{r=2m+1}^{\infty}\alpha_r\phi^r\right),$$
comparing coefficients on both sides gives $a_r=\alpha_r$ for $0\leq r\leq 2m$. Here we used that
$2m=2\left[c/24\right]+2\leq \left[c/8\right] = k $ for $c\geq 32$. 
The coefficients of $\chi_{1/2}^{24m-2c-1}$ are nonnegative because  
$24\cdot 2m-2c-1=48\left[c/24\right]-2c+47\geq 0$. It follows from Lemma~\ref{coeffpos} below that all
coefficients of $\chi'_{M_c}\chi_{1/2}-2c\cdot \chi_{M_c}\chi'_{1/2}$ are negative for $c\geq 1.01$.
Equation~(\ref{ar}) gives now $a_{2m}<0$.

\smallskip
For the second estimate of $a_{2m}$, we obtain from (\ref{schattenentw}) the equation
\begin{equation}\label{spitze1exp1}
 \sum_{r=0}^ka_r(\widetilde{\chi}_{1/2}^{24})^{k-r}=
      q^{c/12-[c/8]} \widetilde{\chi}_{1/2}^{24k-2c} \sum_{n=0}^{\infty} B_n\,q^n. 
\end{equation}
Let 
\begin{equation}\label{spitze1exp2}    
  q^{n+c/12-[c/8]} \widetilde{\chi}_{1/2}^{24k-2c}=
  \sum_{r=0}^{\infty}\beta_{n,r}\, \widetilde{\phi}^r 
\end{equation}
be the expansion of $  q^{c/12-[c/8]} \widetilde{\chi}_{1/2}^{24k-2c}\,q^n$ in powers
of $\widetilde{\phi}=\widetilde{\chi}_{1/2}^{24}=2^{12}q+O(q^2)$.
Using again the B\"urmann-Lagrange theorem, we have  that  $\beta_{n,r}$  is for $r>0$ the coefficient of $q^{r-1}$ in
$$\frac{1}{r}\frac{d( \widetilde{\chi}_{1/2}^{24k-2c}q^{n+c/12-[c/8]})}{dq}
 \left(\frac{q}{\widetilde{\phi}}\right)^r\qquad\qquad\qquad\qquad\qquad\qquad\qquad\qquad\qquad\qquad$$
$$\qquad\qquad=\frac{1}{r}q^{n+c/12-k-1+r}\widetilde{\chi}_{1/2}^{24k-2c-1-24r}
     \left[ (24k-2c)q \widetilde{\chi}'_{1/2}+(n+\frac{c}{12}-k)\widetilde{\chi}_{1/2}\right]$$
and $\beta_{0,0}=2^{24[c/24]-c+24}$.
The coefficients of $\widetilde{\chi}_{1/2}^{24k-2c-1-24r}$ for $r=k-2m$ are nonnegative because  
$24k-2c-1-24r=48\left[c/24\right]-2c+47\geq 0$. With $\widetilde{\chi}_{1/2}$ also 
$(24k-2c)q \widetilde{\chi}'_{1/2}+(n+\frac{c}{12}-k)\widetilde{\chi}_{1/2}$ has for all $n\geq 0$ positive
coefficients. Thus $\beta_{n,k-2m}$ is positive.

Comparing equation~(\ref{spitze1exp1}) and~(\ref{spitze1exp2}) gives 
$a_{k-r}=\sum_{n=0}^{r}\beta_{n,r}B_n$ and hence $a_{2m}\geq 0$,
the desired contradiction. \eop

\smallskip

\noindent{\bf Remark:} The analogous result for unimodular lattices was obtained
by Rains and Sloane~\cite{RaSl-lattice}. 
For even self-dual binary codes and even Kleinian codes the corresponding results can 
be found in~\cite{Ra-shadow}.

\smallskip

\begin{lemma}\label{coeffpos}
The coefficients of the series $2c\chi_{M_c}\chi'_{1/2}-\chi'_{M_c}\chi_{1/2}$
are all positive for $c\geq 1.01$.
\end{lemma}
\pf Let $A=p^{1/24}\chi_{1/2}=\prod_{n=0}^\infty (1+p^{2n+1})$ and $B=p^{c/12}\chi_{M_c}=\prod_{n=2}^\infty (1-p^{2n})^{-1}$.
Then 
$$2c\chi_{M_c}\chi'_{1/2}-\chi'_{M_c}\chi_{1/2}=p^{-c/12-1/24}(2cBA'-B'A)=p^{-c/12-1/24}B(2cA'-A\cdot B'/B)$$ 
and 
it is enough to show that $2cA'-A\cdot B'/B$ has positive coefficients. One has
\begin{eqnarray}\nonumber
B'/B & = & \frac{d}{dp}\left(\sum_{n=2}^{\infty}-\log(1-p^{2n})\right) \\ \nonumber
 & = & \sum_{n=2}^{\infty}\frac{2n\, p^{2n-1}}{1-p^{2n}} \\ \nonumber
 & = &  \sum_{n=2}^{\infty}2(\sigma(n)-1)p^{2n-1},
\end{eqnarray}
where $\sigma(n)$ is the sum of the positive divisors of $n$. Thus the coefficients of $B'/B$
can be estimated from above by the coefficients of $2\sum_{n=1}^{\infty}\sigma(n)p^{2n-1}$.
Let $A=\sum_{n=0}^{\infty}f(n)p^n$. Then the coefficient of $p^n$ in 
$AB'/B$ is estimated from above by $2$ times
\begin{eqnarray}\label{estimation}\nonumber
\sum_{k=1}^{[\frac{n+1}{2}]}\sigma(k)f(n-(2k-1)) & = & 
     \sum_{k=1}^{[\frac{n+1}{2}]}f(n+1-2k)\sum_{d|k}\frac{k}{d} \\ \nonumber
& = & \sum_{d=1}^{[\frac{n+1}{2}]}\frac{1}{d}\sum_{{\scriptstyle k\leq [\frac{n+1}{2}]}\atop {\scriptstyle  d|k}} f(n+1-2k)k \\ \nonumber
& = & \sum_{d=1}^{[\frac{n+1}{2}]}\sum_{r=1}^{[[\frac{n+1}{2}]/d]} f(n+1-2rd)\cdot r \\ \nonumber
& \sim &  \sum_{d=1}^{[\frac{n+1}{2}]} \int_0^{n/2d}f(n-2td)\cdot t \cdot dt  \\ \nonumber
& = &  \sum_{d=1}^{[\frac{n+1}{2}]}\frac{1}{d^2} \int_0^{n/2}f(n-2x)\cdot x \cdot dx \\ 
& < & \frac{\pi^2}{6} \cdot  \int_0^{n/2}f(n-2x)\cdot x \cdot dx.
\end{eqnarray}
The coefficient $f(n)$ counts the number of partitions of $n$ into odd and unequal parts
and one has for $n\longrightarrow \infty$ the asymptotic formula (\cite{Hagis}, Corollary of Thm.~6)
$$f(n)=\sqrt{6}\,(24n-1)^{-3/4}\exp\left(\pi\sqrt{24n-1}/12\right)\left(1+{\rm O}(n^{-1/2})\right). $$
Using this approximation for $f(x)$, the integral  $\int_0^{n/2}f(n-2x)\cdot x \cdot dx$ can be evaluated
explicitly and one obtains $\lim_{n\longrightarrow\infty}\left(\int_0^{n/2}f(n-2x)\cdot x \cdot dx\right)/(n f(n))=6/\pi^2$.
Hence the coefficients of $AB'/B$ are smaller than the coefficients $(n+1)f(n+1)$
in $2cA'$ for large $n$.

It is now also straightforward to justify the approximation of the sum by the integral in~(\ref{estimation}):
The function $f(n-2t)\cdot t$ is not monotone on $[1,n/2]$ but has a single maximum at $t_0\sim \sqrt{24}/(2\pi)\sqrt{n}$.
The possible approximation error is therefore not larger then $\sum_{d=1}^{(n+1)/2}\frac{1}{d}\cdot f(n-2t_0) t_0\sim
\log(n/2)f(n-2t_0) t_0$. But 
$$\lim_{n\longrightarrow\infty}\log(n/2)f(n-2t_0)t_0/(nf(n))=0.$$

We skip the explicit computation of an $N_0$ such that the Lemma holds for $n\geq N_0$.
For $n<3000$ we checked the Lemma directly.
\eop

\bigskip

For smaller values of the central charge~$c$, one can often improve the general upper
bound of~Theorem~\ref{newbound}.  Table~\ref{minsvoas} lists our results for $c\leq 48$.

\smallskip

We use that the character of a \svoa and of its
shadow must have nonnegative integral coefficients.
(More precisely, the dimensions of the Virasoro primaries must be nonnegative.)

As an example, assume that a \svoa of central charge $c=33\ha$ and minimal weight~$2$
exists. This would imply $a_0=1$, $a_1=-67$, $a_2= 670$, $a_3=-201$ and hence
\begin{eqnarray}\nonumber
\chi_V &= & q^{-33\ha/24}\bigl(1+(56816 + a_4)\,q^2+(2072444 - 29 a_4)\,q^{5/2}+\cdots\bigr),\\ \nonumber
\chi_{V'} &= & q^{-33\ha/24}\bigl(a_4/32768\,    q^{3/16} + (823296 - 29 a_4)/32768\,q^{19/16}+\cdots\bigr).
\end{eqnarray}
The initial term of $\chi_{V'}$ gives $a_4\geq 0$ and $2^{15}\mid a_4$. The second coefficient of $\chi_{V'}$
is only nonnegative for $a_4=0$ in which case it is nonintegral.

This kind of arguments gives the upper bounds for all central charges $c$ listed besides  
for $c=10$, $11$, $12\ha$, $13$, $13\ha$, $14\ha$
and $16\ha$, for which we used that the list of \svoas given in Thm.~5.3.2~\cite{Ho-dr}
is complete, which in turn depends on Schellekens classification~\cite{schellekens1} 
of self-dual \voas of central charge~$24$.  However, the method of Schellekens can also 
directly be applied for $c\leq 16\ha$ to show that the obtained value for $\dim V_1$ 
cannot be realized.

\bigskip

It is known that any unimodular lattice in dimensions 
$n\equiv 0\!\pmod{24}$ meeting the upper bound analogous to Theorem~\ref{newbound} must 
be even (see~\cite{Gauter-even}). A similar result holds 
for even self-dual binary codes (cf.~\cite{RaSl-lattice}). For \svoas
of small central charges $c\equiv 0\!\pmod{24}$ we can show:

\begin{theorem}\label{maxodd}
A \svoa $V$ of central charge $c=24$, $48$, $72$ or $96$ with minimal weight
$\mu(V)\geq c/24 + 1/2$ must be a \voa.
\end{theorem} 
\pf For a \svoa which is not a \voa, the first coefficient  of the character
$\chi_{V'}=  q^{-c/24}\sum_{n=0}^{\infty} B_n\,q^n$ vanishes since the conformal
weight of the shadow  is  positive. This is impossible for $c=24$. For the other
values of~$c$ it is impossible to find values of the still undetermined $a_i$ such that
all $B_n$, $n\geq 1$,  are nonnegative.
\eop

The largest minimal weight of a \svoa $V$ which is not a \voa for the above four values of $c$
is therefore $c/24 + 1/2$.

\smallskip

\noindent{\bf Problem:}  Is a self-dual \svoa $V$ of central charge
\hbox{$c\equiv 0\!\!\! \pmod{24}$} with minimal weight $\mu(V)= \left[\frac{c}{24}\right] + 1 $ always a
vertex operator algebra? (The proof in~\cite{Gauter-even} for lattices cannot directly be generalized
since it uses a lower bound for the minimal norm of the shadow of an odd unimodular lattice
whereas the analogous bound for \svoas is not obvious.) 

\medskip

The case of self-dual \svoas $V$ of central charge~$48$ will be investigated
further in Section~\ref{c48}.

\begin{table}\caption{Highest minimal weight $\mu$ of a self-dual \svoa of central charge~$c$ for $c\leq 48$}\label{minsvoas}
{\small
$$\begin{array}{l||c|c|c|c|c|c|c|c|c|c|c|c|c|c|c|c}
c   & \ha & 1   & 1\ha & 2   & 2\ha & 3   & 3\ha & 4   & 4\ha & 5   &  5\ha & 6   & 6\ha & 7   & 7\ha & 8 \\\hline
\mu & \ha & \ha & \ha  & \ha &  \ha & \ha &  \ha & \ha &  \ha & \ha &  \ha  & \ha &  \ha & \ha &  \ha & 1 \\ 
\\
c   & 8\ha & 9   & 9\ha & 10  & 10\ha & 11  & 11\ha & 12  & 12\ha & 13  & 13\ha  & 14  & 14\ha & 15  & 15\ha & 16 \\\hline
\mu & \ha  & \ha & \ha  & \ha & \ha   & \ha &  \ha  & 1   &  \ha  & \ha &  \ha   & 1   &  \ha  & 1   &   1   & 1  \\ 
\\
c   & 16\ha & 17   & 17\ha & 18  & 18\ha & 19  & 19\ha & 20  & 20\ha & 21  & 21\ha  & 22  & 22\ha & 23  & 23\ha   & 24 \\\hline
\mu & \ha   & 1    & 1     & 1   & 1     & 1   &  1    & 1   &  1    & 1   & 1      & 1   &  1    & 1   &  1\ha   & 2  \\
\\
c   & 24\ha    & 25     & 25\ha & 26     & 26\ha &  27   & 27\ha &  28    & 28\ha & 29    & 29\ha  & 30      & 30\ha   & 31      & 31\ha    & 32 \\\hline
\mu &\ha\bis1\ha &1\bis1\ha&1\bis1\ha&1\bis1\ha&1\bis1\ha&1\bis1\ha&1\bis1\ha&1\bis1\ha&1\bis1\ha&1\bis1\ha&1\bis1\ha  &1\bis1\ha&1\bis1\ha &1\bis1\ha &1\ha & 2  \\ 
\\
c   & 32\ha    & 33      & 33\ha     & 34     & 34\ha   & 35     & 35\ha   & 36       & 36\ha   & 37   & 37\ha    & 38      & 38\ha   & 39      & 39\ha    & 40 \\\hline
\mu & 1\bis1\ha&1\bis1\ha&1\bis 1\ha &1\bis 2 &1\bis 2  &1\bis 2 &1\bis 2  &1\bis 2   &1\bis 2  &1\bis 2 &1\bis 2 &1\bis 2  &1\bis 2  &1\bis 2  &1\ha\bis 2   & 2  \\ 
\\
c   & 40\ha  & 41      & 41\ha   & 42     & 42\ha   & 43     & 43\ha   & 44       & 44\ha    & 45     & 45\ha   & 46      & 46\ha   & 47        & 47\ha     & 48 \\\hline
\mu & 1\bis 2 &1\bis 2  &1\bis 2  &1\bis 2 &1\bis 2  &1\bis 2 &1 \bis 2 &1\bis 2   &1 \bis 2  &1\bis 2 &1\bis 2  &1\bis 2  &1\bis 2  & 1\ha\bis 2&1\ha\bis 2 & 2\bis 3  \\ 
\end{array}$$}
\end{table}

\bigskip

Examples of \svoas achieving the lower bounds given in Table~\ref{minsvoas} can be constructed as follows:

For $c<8$ and all other $c$ with $\mu=\halb$  one uses $V_{\rm Fermi}^{\otimes 2c}$ where
$V_{\rm Fermi}$ is the $c=\halb$ \svoa describing a single fermion.

For $8 \leq  c \leq 23$ and $\mu=1$ we use Thm.~5.3.2~\cite{Ho-dr} and the table
in Section~2 of~\cite{Ho-shadow}. Both results
depend on Schellekens classification of self-dual $c=24$ \voas and
use the construction of Chapter~3 in~\cite{Ho-dr}. In particular, for integral $c$ these examples
are just the lattice \svoas associated to an integral unimodular lattice with 
minimal norm $2$ or $3$.

For $c=23\halb$ the shorter Moonshine module $V\!B^{\natural}$ of~\cite{Ho-dr}, Ch.~4, 
is an example with $\mu=1\halb$ and for $c=24$ the Moonshine module $V^{\natural}$~\cite{FLM} 
is an example with $\mu=2$.

For $c=32$, $40$ and $48$ the ${\bf Z}_2$-orbifolds of lattice \voas associated to extremal even unimodular
lattices in dimension $32$, $40$ and $48$~\cite{DGM} are examples with $\mu=2$. Since for $c=32$ and
$40$ these \voas have a Virasoro frame (cf.~\cite{DGH-virs}), one uses the same construction as for $V\!B^{\natural}$
to obtain self-dual \svoas of central charge~$31\ha$ and~$39\ha$ with $\mu=1\ha$. 

All remaining examples can be obtained by taking tensor products of two \svoas $V$ and $W$ 
of smaller central charge and using $\mu(V\otimes W)=\min\{\mu(V),\,\mu(W),\,2\}$.


\section{Self-dual \N1 supersymmetric \svoas{} }\label{N1svoas}

We recall that an \N1 supersymmetric \svoa $V$ is an \svoa $V$ together with a superconformal element 
$\tau\in V_{3/2}$ such that the operators $G_{n+1/2}=\tau_{n+1}$ generate a representation of
the Neveu-Schwarz superalgebra on $V$. This is the case precisely if $\tau_2 \tau=\frac{2}{3} c{\bf 1}$,
$\tau_1 \tau=0$ and $\tau_0 \tau=2\omega$. The $\sigma$-twisted modules 
(or Ramond sectors) of $V$ admit then a representation
of the Ramond superalgebra.
For a $\sigma$-twisted module one has under some unitary assumption that the conformal weight satisfies 
$h\geq \frac{c}{24}$ (cf.~\cite{LuTh}, p.~242).	  

\bigskip

\begin{definition} {\rm The {\it minimal superconformal weight\/} $\mu^*(V)$ of an \N1
\svoa $V$ is defined as the smallest degree of a highest weight vector of the
\N1 super Virasoro algebra different from the vacuum vector. In case the only 
highest weight vector is the vacuum vector, we let $\mu^*(V)=\infty$.}
\end{definition}

One has
$$\chi_V=q^{-c/24}\Biggl(\prod_{n=2}^{\infty}\frac{1+q^{n-1/2}}{1-q^n}+
\prod_{n=1}^{\infty}\frac{1+q^{n-1/2}}{1-q^n}\biggl(\sum_{\ i\geq \mu^*(V)\ } P_i\cdot q^i
\biggr) \biggr),$$
where $P_i$ is the dimension of the space of highest weight vectors of degree~$i$ for the
\N1 super Virasoro algebra.

\medskip
For a self-dual \N1 supersymmetric \svoas $V$, the shadow $V'$ is the unique $\sigma$-twisted $V$-module and hence has
conformal weight $h(V')\geq c/24$. This implies that $\chi_{V'}$ has no pole in the cusp $i\infty$.
By using the relation between the characters of $V$ and $V'$ given in the introduction 
it follows that
$\chi_V$ has no pole in the cusp~$1$. Hence equation~(\ref{svoachar}) gives
$$ 
\chi_V=\sum_{r=0}^ka_r\chi_{1/2}^{2c-24r}=q^{-c/24}\sum_{n=0}^{\infty} C_{n}\,q^{n/2}
$$ 
with $k=\left[\frac{c}{12}\right]$, function $\chi_{1/2}$ as in~(\ref{chi12}), 
and uniquely determined integers $a_0$, $\ldots$, $a_k$.

Denote by $\chi_{M_c^{N\!=\!1}}=q^{-c/24}(\prod_{n=2}^{\infty}\frac{1+q^{n-1/2}}{1-q^n})$ the character
of the \N1 Viraso vertex operator superalgebra generated by $\tau$.
Following~\cite{Witten-3dgravity} we make the following definition:
\begin{definition}\label{n1extremal}{\rm 
If the $a_0$, $\ldots$, $a_k$ are chosen such that one has
\begin{equation}\label{n1extchar}
\chi_V=\chi_{M_c^{N\!=\!1}}\cdot(1+A_{2k+1}\,q^{(k+1)/2}+A_{k+2}\,q^{(k+2)/2}+\cdots)
\end{equation}
then~(\ref{n1extchar}) is called the {\it extremal (\N1 supersymmetric) character\/} 
for central charge~$c$. A self-dual \N1 supersymmetric \svoa with 
this character is called an {\it extremal \N1 supersymmetric \svoa.}}
\end{definition}
It follows directly from the definitions that the minimal superconformal weight of a
self-dual \N1 supersymmetric \svoa $V$ satisfies $ \mu^*(V)\geq \halb\left[\frac{c}{12}\right]+\halb.$
 
The character of the shadow $V'$ of a self-dual \N1 supersymmetric \svoa  $V$ is 
$$ 
\chi_{V'}=\alpha\, \sum_{r=0}^ka_r\widetilde{\chi}_{1/2}^{2c-24r}=
    \alpha\,  q^{c/12-[c/12]}\sum_{n=0}^{\infty} B_n\,q^n
$$ 
with $\alpha=1$ for integral  and $\alpha=1/\sqrt{2}$ for nonintegral values of $c$
and  $\widetilde{\chi}_{1/2}$ as in~(\ref{chis12}). If the $a_i$, $0\leq i \leq k$, are
chosen as in Definition~\ref{n1extremal} we call $\chi_{V'}$ the {\it extremal shadow character.}

\begin{theorem}\label{n1upperbound}
A self-dual \N1 supersymmetric \svoa has minimal superconformal weight 
$$ \mu^*(V)\leq \halb\left[\frac{c}{12}\right]+\halb$$
unless $c=23\ha$ in which case $\mu^*(V)\leq 3/2$.
\end{theorem}
\pf Let $k=[c/12]$ and $p=q^{1/2}$. We compute the coefficient $A_{k+1}$ of the extremal character 
$\chi_V=\chi_{M_c^{N\!=\!1}}\cdot(1+A_{k+1}\,q^{(k+1)/2}+\cdots)$ by the method used
in the proof of Theorem~\ref{newbound} (see~(\ref{ar}) and the following equation) 
and obtain that $A_{k+1}$ is the coefficient of $p^k$ in
\begin{equation}\label{n1ar}
-\frac{1}{k+1}\,p^{k+1}\cdot\chi_{1/2}^{24(k+1)-2c-1}
  \left[\chi'_{M_c^{N\!=\!1}}\chi_{1/2}-2c\cdot \chi_{M_c^{N\!=\!1}}\chi'_{1/2}\right].
\end{equation}
The coefficients of $p^{24(k+1)-2c}\chi_{1/2}^{24(k+1)-2c-1}=1+O(p)$ are nonnegative since 
$24(k+1)-2c-1= 23-2c+24[c/12]\geq 0$. Similar as in Lemma~\ref{coeffpos} it can be shown that all coefficients 
of $-p^{(2c+1)/24}(\chi'_{M_c^{N\!=\!1}}\chi_{1/2}-2c\cdot \chi_{M_c^{N\!=\!1}}\chi'_{1/2})$ are
positive besides the coefficient of $p$ which is zero.

Thus $A_{k+1}$ is positive unless $c=23\ha$ in which case the extremal character
is $q^{-47/48}(1+4372\,q^{3/2}+O(q))$ yielding $\mu^*\leq 3/2$. \eop

\smallskip

\noindent{\bf Remark:} Witten defines extremal \N1 supersymmetric \svoas 
by only requiring that $\mu^*\geq \halb[c/12]$. 
Then the character is only determined up to the addition of a constant. 
The reason for this modification is that for $c$ a multiple of $24$ and $c\geq 96$
the first coefficient of the expansion of the extremal shadow character becomes negative, i.e., 
for those values of $c$ an extremal \N1 supersymmetric \svoa cannot exist. Hence for those $c$
the bound in Theorem~\ref{n1upperbound} can be improved by~$\halb$.

Numerical evidence suggests
that besides for the mentioned case all coefficients of the extremal character and
of the extremal shadow character are positive integral numbers.

\bigskip

Examples of extremal self-dual \N1 supersymmetric \svoas are
$V_{\rm fermi}^{\otimes 2c}$ for $c=(3/2)k$, $k=1$, $\ldots$, $7$, the \svoa
$V_{D_{12}^+}$~\cite{Duncan-conway} and the odd Moonshine module~$V\!O^{\natural}$~\cite{DGH-Monster}.

There exist other self-dual \svoas with the right extremal character who might have also
the additional structure of an \N1 supersymmetric \svoa.


\section{Central Charge~$c=48$ and the Monster}\label{c48}

In this section, we study self-dual \voas and \svoas of central charge~$48$ 
with large minimal conformal or superconformal weight. The question of possible
monster symmetries of such \svoas is investigated.

\medskip

We recall that Theorem~\ref{maxodd} implies for central charge~$48$  
that the minimal weight of a \svoa $V$ which is not a \voa can be at most $5/2$.
More precisely, we can deduce from the requirement that the character of $V$ and 
its shadow $V'$ have nonnegative integral coefficients that for 
a \svoa with minimal weight $5/2$ there are only the following two possibilities 
for the characters:
\begin{eqnarray}\nonumber
\chi_{V} & = & q^{-2}+ 1 +  192512\, q^{1/2} +  21590016\, q+ 863059968\, q^{3/2} +  20256751892\, q^{2} + \cdots\\ \label{c48short}
\chi_{V'} & = & q^{-1}+ 1 + 42991892\, q+  40491808768\, q^{2}+ 8504047840194\, q^{3}+\cdots
\end{eqnarray}
\begin{eqnarray}\nonumber
\chi_{V} & = & q^{-2}+ 1 +  196608\, q^{1/2} + 21491712\, q+ 864288768\, q^{3/2} + 20246003988\, q^{2} 
+\cdots\\\label{c48long}
\chi_{V'} & = & 25+ 42991616\, q+ 40491810816\, q^{2}+ 8504047828992\, q^{3}+\cdots
\end{eqnarray}

For lattices of dimension~$48$ and codes of length~$48$ an analogous result can be found
in~\cite{HKMV-48}. There exist lattices as well as codes realizing both possibilities
for the theta series respectively weight enumerator.
However, for \voas the first case can be excluded.
\begin{theorem}\label{noneighbour}
Let $V$ be a \svoa of central charge~$48$ with minimal weight $\mu(V)=5/2$. 
Then the shadow of $V$ has minimal conformal weight $2$. 
\end{theorem}
\pf Since $\mu(V)=5/2$ is nonintegral, $V$ is not a \voa.
Hence the even part $V_{(0)}$ has four irreducible modules $V_{(0)}$, $V_{(1)}$, $V_{(2)}$ and $V_{(3)}$,
where $V=V_{(0)}\oplus V_{(1)}$ and $V'=V_{(2)}\oplus V_{(3)}$. The simple current
extensions $W=V_{(0)}\oplus V_{(2)}$ and $\widetilde{W}=V_{(0)}\oplus V_{(3)}$ can be given the structure 
of self-dual \voas and both are neighbours of $V$. 
Assume that the shadow of $V$ has minimal conformal weight different from $2$.
Then the character of $V'=V_{(2)}\oplus V_{(3)} $ is given by equation~(\ref{c48short}). 
Thus either $V_{(2)}$ or $V_{(3)}$ has a one-dimensional degree~$1$ part; say $V_{(3)}$. 
In this case $\widetilde{W}_1$ is one-dimensional and generates a Heisenberg \subvoa of central charge~$1$ with graded
dimension $\prod_{n=1}^{\infty}(1-q^n)^{-1}$. This implies $\dim \widetilde{W}_2\geq 2$ 
and hence $\dim \widetilde{W}_2= 2$
using~(\ref{c48short}) again. A Heisenberg \voa contains the  one-parameter family 
$\omega_{\lambda}=\halb h_{(-1)}^2{\bf 1}+ \lambda h_{(-2)}{\bf 1}$,  $\lambda\in \C$,
of possible Virasoro vectors, where $\omega_{\lambda}$ has central charge $1-12\lambda^2$. 
For the Virasoro element $\omega_{\lambda_0}$ of $\widetilde{W}$ one has therefore $\lambda_0\not=0$ 
since $\widetilde{W}$ has central charge~$48$. 
The \voa $\widetilde{W}=V_{(0)}\oplus V_{(3)}$ admits a natural automorphism $\sigma$ which acts 
by multiplication  with $+1$ on $V_{(0)}$ and with $-1$ on $V_{(3)}$. 
Since $\widetilde{W}_1=\C h_{(-1)}{\bf 1}\subset V_{(3)}$, it follows that  $h_{(-1)}^2{\bf 1}$
is in the $+1$-eigenspace of the involution $\sigma$. But $\omega_{\lambda_0}$ is also in the $+1$-eigenspace.
This contradicts $\dim (V_{(0)})_2=1$. 

Our assumption about the minimal conformal weight of the shadow to be different from $2$ 
must therefore be wrong.
\phantom{xx}  \eop

If a \svoa with minimal weight $5/2$ and minimal conformal weight of the shadow be $1$
would exist, then the neighbour \voa $W$ as in the proof would be an example of an 
extremal \voa, i.e., a self-dual \voas with minimal weight $\mu=3$.

\medskip

It follows from Theorem~\ref{n1upperbound} or the discussion
in~\cite{Witten-3dgravity}, Sect.~3.2, that the minimal superconformal weight
of an \N1 supersymmetric \svoa $V$ of central charge~$48$,
is at most $5/2$. If $V$ is extremal, 
the characters of $V$ and $V'$ are (\cite{Witten-3dgravity}, eq.~(3.60)):
\begin{eqnarray}\nonumber
\chi_{V}  & = & q^{-2}+ q^{-1/2}+ 1+ 196884\,q^{1/2}+ 21493760\,q+ 864299970\,q^{3/2} \\ \nonumber
& & \qquad\qquad  +  20246053140 \,q^{2} + 333202640600 \,q^{5/2}  + 4252023300096 \,q^{3}  +\cdots,\\ \label{extc48N1}
\chi_{V'} & = & 1 +  42987520\,q+ 40491712512\,q^{2}+ 8504046600192\,q^{3} +\cdots. 
\end{eqnarray}
\begin{theorem}[cf.~\cite{Witten-3dgravity}, discussion at the end of Sect.~3.3]\label{extN1neighbour}
Let $V$ be an extremal self-dual \N1 supersymmetric \svoa of central charge~$48$. 
Then $V$ has an extremal self-dual \voa $W$ as neighbour. 
\end{theorem} 
\pf
As in the proof of Theorem~\ref{noneighbour}, let $W=V_{(0)}\oplus V_{(2)}$ and 
$\widetilde{W}=V_{(0)}\oplus V_{(3)}$ be the two \voa neighbours of $V$.
From the character of $V'=V_{(2)}\oplus V_{(3)}$ given in equation~(\ref{extc48N1}) it
follows that either $V_{(2)}$ or $V_{(3)}$, say $V_{(2)}$, has zero-dimensional degree~$2$ part.
Then $V_{(2)}$ has minimal conformal weight at least~$3$ and $W=V_{(0)}\oplus V_{(2)}$
has minimal weight~$3$, i.e., is an extremal \voa of central charge~$48$.
\phantom{xx} \eop

\noindent{\bf Remark:} The same argument shows that a neighbour of an extremal self-dual \N1 
supersymmetric \svoa of central charge~$72$ is an extremal \voa. The characters of 
$V$ and $V'$ of such a \svoa can be found in~\cite{Witten-3dgravity}, Appendix~A.

\smallskip

Witten observes that the first coefficients of the modular functions in~(\ref{extc48N1})
are simple linear combinations of dimensions of irreducible representations
of the monster simple group and further that such a decomposition is compatible
with the  \N1 super Virasoro algebra module structure of $V$ (see~\cite{Witten-3dgravity},
eq.~(3.61) and (3.62)).

Since any action of the monster group on a central charge~$48$ \svoa $V$ induces 
such an action on the $V_{(0)}$-modules $V_{(i)}$, $i=0$, $1$, $2$, $3$, and 
hence on the \voa neighbours $W$ and  $\widetilde{W}$ (it is easy to see that 
this action respects the \voa structure), we have the following corollary 
to Theorem~\ref{extN1neighbour}:

\begin{cor}\label{monsterN1}
Assume that $V$ is an extremal  self-dual \N1 \svoa of central charge~$48$ 
on which the monster acts by automorphisms. 
Then the monster acts also on the extremal neighbour \voa $W$ as in Theorem~\ref{extN1neighbour} by automorphisms
such that the actions  coincides on the common \subvoa  $V_{(0)}$ of $V$ and $W$.
\end{cor}

\medskip

Witten asks in~\cite{Witten-3dgravity}, Sect.~3.1, if extremal self-dual \voas
of central charge $c=24k$ exist for all natural numbers $k$, if they are unique,
and if they have a monster symmetry. In the following, it is shown   
that for $c=48$ at least a monster symmetry is impossible under certain 
natural assumptions. 

\smallskip

The character of a self-dual \voa of central charge divisible by~$24$ is a polynomial with integer coefficients
in the modular invariant $j=q^{-1}+ 744+ 196884\,q+21493760\,q^2+\cdots$ or, equivalently, in
$J=j-744$, the character of the moonshine module $V^{\natural}$. In particular, 
for the character of an extremal \voa $W$ of central charge~$48$ one has (see~\cite{Ho-dr}, Table~5.1):
\begin{eqnarray}\label{extc48}
\chi_W & = & J^2 - 393767\\\nonumber
 & = & q^{-2}  + 1 + 42987520\, q + 40491909396\, q^2  + 8504046600192\, q^3+\cdots . 
\end{eqnarray}
If we assume that the monster acts on $W$ nontrivially by \voa automorphisms,
the simplest way to consider $W$ as a module for the monster
is to assume that one has 
\begin{equation}\label{monsterrep}
W=V^{\natural}\otimes V^{\natural} - (2R_2 + R_1) 
\end{equation} 
as graded monster modules, where $R_1$ denotes the trivial $1$-dimensional and 
$R_2$ is the irreducible $196883$-dimensional representation of the monster.
This also guarantees that the monster module structure is compatible
with the Virasoro module structure of $W$ as one can easily check. 
For an element $g$ in the monster the graded trace of $g$ acting on $W$ is then given by
\begin{equation}\label{traceid}
 \tr(g|W)=q^{-2}\sum_{n=0}^{\infty} \tr(g|W_n)\,q^n=T_g^2-(2\,\tr(g|R_2) + 1),
 \end{equation}
where $T_g$ is the McKay-Thompson series of~$g$, i.e., the graded trace of $g$ acting 
on~$V^{\natural}$.

Furthermore, if we assume that the monster module structure of the first few homogeneous spaces
of a \N1 supersymmetric \svoa $V$ of minimal superconformal weight $5/2$ is 
the one given in~\cite{Witten-3dgravity}, eq.~(3.61) and (3.62), then the 
monster module structure of the extremal neighbour of $V$ as in Corollary~\ref{monsterN1}
is also compatible with~(\ref{monsterrep}) at least if we modify
the monster module structure of~\cite{Witten-3dgravity}~eq.~(3.61) and (3.62) by 
exchanging $R_1+R_4+R_5$ with $R_3+R_6$ if necessary. 
(It was remarked in~\cite{Witten-3dgravity} that such modifications are possible.)
More precisely, we could assume that for $V$ one has
\begin{eqnarray} \nonumber
V_0     & = &  R_1  \\ \nonumber
V_{1/2} & = &  0  \\ \nonumber
V_1     & = &  0  \\ \nonumber
V_{3/2} & = &  R_1  \\ \nonumber
V_2     & = &  R_1  \\ \nonumber
V_{5/2} & = &  R_1+R_2  \\ \nonumber
V_3     & = &  R_1+R_2+R_3  \\ \nonumber
V_{7/2} & = &  2R_1+2R_2+R_3+R_4  \\ \nonumber
V_4     & = &  4R_1+4R_2+R_3+2R_4+R_5   \\ \nonumber
V_{7/2} & = &  5R_1+5R_2+2R_3+3R_4+2R_5+R_7  \\ \nonumber
V_5     & = &  5R_1+7R_2+4R_3+4R_4+2R_5+2R_6+R_7+R_8 \qquad\qquad
\end{eqnarray}
and for $V'$ one has
\begin{eqnarray}\nonumber
V'_0     & = & 0  \\ \nonumber
V'_1     & = & 0  \\ \nonumber
V'_2     & = &  R_1  \\ \nonumber
V'_3     & = & 2\times( R_1+R_2+R_3)  \\ \nonumber
V'_4     & = & 2\times( 2R_1 + 3R_2 + 2R_3 + R_4 + R_6)\\ \nonumber
V'_5     & = & 2\times( 3R_1 + 7R_2 + 6R_3 + 2R_4 + 4R_6 + R_7 + R_8),\qquad\qquad
\end{eqnarray}
where $R_i$ denotes the $i$-th representation of the monster.
This decomposition remains compatible with an \N1 super Virasoro algebra module structure.

We will also use the following conjecture about the structure of ${\bf Z}_2$-orbifolds of 
self-dual \voas{} which seems not to be completely proven: 
\begin{conjecture}\label{orbifold}{
Let $t$ be an involutive automorphism of a self-dual \voa $W$. Then the fixpoint \subvoa
$W^{\langle t \rangle}$ is rational and has the fusion algebra ${\bf Z}[{\bf Z}_2\times {\bf Z}_2]$. 
The conformal weights of the four isomorphism types of irreducible $W^{\langle t \rangle}$-modules are 
either congruent to $0$, $0$, $0$, $1/2\!\!\pmod{1}$ (case I)  or 
$0$, $0$, $1/4$, $3/4\!\!\pmod{1}$ (case II).}
\end{conjecture}
For $W$ the  self-dual lattice \voa $V_{E_8}$ of central charge~$8$ with
$E_8({\bf C})$ as automorphism group, the two conjugacy classes of involutions
of $E_8({\bf C})$ realize both cases I and II. For $W$ the moonshine module $V^{\natural}$,
the two conjugacy classes of involutions in the monster correspond both to case I.

In case I, one can extend $W^{\langle t \rangle}$ by the module of conformal weight $1/2\!\!\pmod{1}$ to
obtain a self-dual \svoa as neighbour of $W$.

\begin{theorem}\label{nomonsterextvoa}
The monster cannot act by automorphisms on an extremal self-dual \voa $W$ of central charge~$48$
such that as a graded monster module one has  $W=V^{\natural}\otimes V^{\natural} - (2R_2 + R_1) $
provided Conjecture~\ref{orbifold} holds.
\end{theorem}
In fact, we only will need that the monster module structure of $W$ is the 
stated one for $W_n$, $0\leq n \leq 5$.

\smallskip
\pf
Let $t$ be an involution in the monster which has a twofold cover of
the baby monster as centralizer, i.e., an involution of type~$2A$ in  
atlas notation. Then the character of the fixpoint \voa $W^{\langle t \rangle}$ is
\begin{eqnarray}\nonumber
\chi_{W^{\langle t \rangle}} &\!=\! & \halb\left(\chi_W+\tr(t|W)\right) \\ \nonumber
&\! =\! & q^{-2}+1+21590016\, q+20256751892\, q^2+4252454830080\, q^3+ \cdots  
\end{eqnarray}
where we have used equations~(\ref{extc48}), (\ref{traceid}), the character value $\tr(t|R_2)=4371$  
and the $q$-expansion of the Thompson series $T_t$ as conjectured in~\cite{CoNo} and proven 
in~\cite{Bo-lie}.
By comparing this expansion with the one in~(\ref{c48short}) it
follows that $\chi_{W^{\langle t \rangle}}=\chi_{V_{(0)}}$ where $V_{(0)}$ is the even \subvoa
of a \svoa of central charge~$48$ with minimal weight $\mu(V)=5/2$
and shadow of minimal conformal weight~$1$, because both functions are
modular functions for $\Gamma_0(2)$. The width of $\Gamma_0(2)$ in its other cusp is $2$
and hence $\chi_{W^{\langle t \rangle}}(-1/\tau)$ has an expansion in powers of $q^{1/2}$. 
This implies that $t$ is an involution in ${\rm Aut}(W)$ belonging to case I 
of Conjecture~\ref{orbifold}: Denoting the four irreducible modules of $W^{\langle t \rangle}$ by $M_i$,
$i=1$, $\ldots$, $4$,  the expansion of 
$\chi_{W^{\langle t \rangle}}(-1/\tau)=\frac{1}{2}\sum_{i=1}^3\chi_{M_i}(\tau)$
would contain in case II non-even powers of $q^{1/4}$.
Since we are in case I, 
we can extend the \voa $W^{\langle t \rangle}$ to a self-dual \svoa $V$. From the explicit expansion
of $\chi_{W^{\langle t \rangle}}$ in the other cusp, we see that the characters of $V$ and its shadow  are
the ones given in~(\ref{c48short}). Thus $V$ has minimal weight $5/2$ and the shadow of $V$
has minimal weight~$1$. However, by Theorem~\ref{noneighbour} such a \svoa cannot exist.
\eop

\noindent{\bf Remark:}  If we take instead of an involution of type~$2A$ an involution 
of type~$2B$ in the monster, then the corresponding neighbour \svoa has the 
character given in~(\ref{c48long}).

\smallskip

Theorem~\ref{nomonsterextvoa} together with Corollary~\ref{monsterN1} implies
that under reasonable assumptions it is impossible for the monster simple group to act on
an extremal \N1 supersymmetric \svoa in the way suggested in~\cite{Witten-3dgravity}.

If we would have used in Corollary~\ref{monsterN1} for $V$ exactly the monster module
structure as in~\cite{Witten-3dgravity}, then the character of the constructed \svoa neighbour 
of $W$ would have non-integral coefficients.

\bigskip
\noindent{\bf Acknowledgments.} I like to thank M.~Gaberdiel for useful discussions.

\providecommand{\bysame}{\leavevmode\hbox to3em{\hrulefill}\thinspace}
\providecommand{\MR}{\relax\ifhmode\unskip\space\fi MR }
\providecommand{\MRhref}[2]{%
  \href{http://www.ams.org/mathscinet-getitem?mr=#1}{#2}
}
\providecommand{\href}[2]{#2}


\begin{thebibliography}{HKMV05}

\bibitem[AKM07]{greek-ext}
Spyros Avramis, Alex Kehagias, and Constantina Mattheopoulou,
  \emph{{Three-dimensional AdS gravity and extremal CFTs at c=8m}}, preprint
  (2007), arXiv:0708.3386 [hep-th].

\bibitem[Bor92]{Bo-lie}
R.~E. Borcherds, \emph{{Monstrous moonshine and monstrous Lie superalgebras}},
  Invent. math. \textbf{109} (1992), 405--444.

\bibitem[CN79]{CoNo}
J.~H. Conway and S.~P. Norton, \emph{{Monstrous Moonshine}}, Bull. London Math.
  Soc. \textbf{11} (1979), 308--339.

\bibitem[DGH88]{DGH-Monster}
L.~Dixon, P.~Ginsparg, and J.~Harvey, \emph{{Beauty and the Beast:
  Superconformal Symmetry in a Monster Module}}, Commun. Math. Phys.
  \textbf{119} (1988), 221--241.

\bibitem[DGH98]{DGH-virs}
Chonying Dong, Robert Griess, and Gerald H{\"o}hn, \emph{{Framed Vertex
  Operator Algebras, Codes and the Moonshine Module}}, Comm. Math. Phys.
  \textbf{193} (1998), 407--448, q-alg/9707008.

\bibitem[DGM90]{DGM}
L.~Dolan, P.~Goddard, and P.~Montague, \emph{{Conformal Field Theory, Triality
  and the Monster Group}}, Physical Letters B \textbf{236} (1990), 165--172.

\bibitem[DLM98]{DLM-twistedrep}
Chongying Dong, Haisheng Li, and Geoffrey Mason, \emph{Twisted representations
  of vertex operator algebras}, Math. Ann. \textbf{310} (1998), no.~3,
  571--600.

\bibitem[Dun05]{Duncan-conway}
John~F. Duncan, \emph{{Super-moonshine for Conway's largest sporadic group}},
  preprint, to appear in Duke Math. J. (2005), arXiv:math/0502267v3.

\bibitem[FLM88]{FLM}
Igor Frenkel, James Lepowsky, and Arne Meuerman, \emph{{Vertex Operator
  Algebras and the Monster}}, Academic Press, San Diego, 1988.

\bibitem[Gab07]{Gaber-ext}
Mathias Gaberdiel, \emph{{Constraints on extremal self-dual CFTs}}, preprint
  (2007), arXiv:0707.4073 [hep-th].

\bibitem[Gau01]{Gauter-even}
Mark Gaulter, \emph{Minima of odd unimodular lattices in dimension {$24m$}}, J.
  Number Theory \textbf{91} (2001), 81--91.

\bibitem[GY07]{GaiYin-ext}
Davide Gaiotto and Xi~Yin, \emph{{Genus Two Partition Functions of Extremal
  Conformal Field Theories}}, preprint (2007), arXiv:0707.3437 [hep-th].

\bibitem[Hag64]{Hagis}
Peter Hagis, Jr., \emph{Partitions into odd and unequal parts}, Amer. J. Math.
  \textbf{86} (1964), 317--324.

\bibitem[HKMV05]{HKMV-48}
Masaaki Harada, Masaaki Kitazume, Akihiro Munemasa, and Boris Venkov, \emph{On
  some self-dual codes and unimodular lattices in dimension $48$}, European J.
  Combin. \textbf{26} (2005), 543--557.

\bibitem[H{\"o}h95]{Ho-dr}
Gerald H{\"o}hn, \emph{{Selbstduale Vertexoperatorsuperalgebren und das
  Babymonster}}, Ph.D. thesis, {U}niversit{\"a}t {B}onn, 1995, see: {B}onner
  {M}athematische {S}chriften {\bf 286}, arXiv:0706.0236.

\bibitem[H{\"o}h97]{Ho-shadow}
\bysame, \emph{{Self-dual Vertex Operator Superalgebras with Shadows of large
  minimal weight}}, Internat. Math. Res. Notices \textbf{13} (1997), 613--621,
  q-alg/9608023.

\bibitem[H{\"o}h07]{Ho-conformal}
\bysame, \emph{{Conformal Designs based on Vertex Operator Algebras}}, to
  appear in Adv. Math. (2007), arXiv:math/0701626.

\bibitem[LT89]{LuTh}
D.~L\"ust and S.~Theissen, \emph{{Lectures on String Theory}}, Lecture Notes in
  Physics, Springer-Verlag, Berlin, Heidelberg, New York, 1989.

\bibitem[Man07]{Mansch-ext}
Jan Manschot, \emph{{${\rm AdS}_3$ Partition Functions Reconstructed}},
  preprint (2007), arXiv:0707.1159 [hep-th].

\bibitem[Rai98]{Ra-shadow}
E.~M. Rains, \emph{Shadow bounds for self-dual codes}, IEEE Trans. Inform.
  Theory \textbf{44} (1998), 134--139.

\bibitem[RS98]{RaSl-lattice}
E.~M. Rains and N.~J.~A. Sloane, \emph{The shadow theory of modular and
  unimodular lattices}, Journal of Number Theory \textbf{73} (1998), 359--389.

\bibitem[Sch93]{schellekens1}
A.~N. Schellekens, \emph{{Meromorphic $c=24$ Conformal Field Theories}}, Comm.
  Math. Phys. \textbf{153} (1993), 159--185.

\bibitem[Sto99]{Stone-space}
E.C. Stone, \emph{{Communications technologies for space exploration}},
  Proceedings of the IEEE \textbf{87} (1999), 1044--1046.

\bibitem[Wit07]{Witten-3dgravity}
Edward Witten, \emph{Three-dimensional gravity revisited}, arXiv.org:0706.3359
  [hep-th].

\bibitem[Zhu90]{Zhu-dr}
Yongchang Zhu, \emph{{Vertex Operator Algebras, Elliptic Functions, and Modular
  Forms}}, Ph.D. thesis, Yale University, 1990, {appeared as: {\it Modular
  invariance of characters of vertex operator algebras,\/} J. Amer. Math. Soc
  {\bf 9} (1996)}.

\end{thebibliography}
\end{document}